\documentclass[12pt]{article}
\usepackage{amsmath}
\usepackage{amssymb}
\usepackage{geometry} 
\title{A new criterion of the Riemann hypothesis}
\author{Roupam Ghosh}
\date{\today} 

\begin{document}

\maketitle
\begin{center}
\textit{Abstract:}
\\
We introduce a new criterion which if satisfied implies the Riemann hypothesis.\end{center}
\section{Background}
The Riemann Hypothesis is usually stated as, the non-trivial zeros of the Riemann zeta function lie on the line $\Re(s) = \frac{1}{2}$. Although, this is the standard formulation, one of the exciting features of this problem is, it can be formulated in many different and unrelated ways. One of them was derived by Nyman and Beurling, and now known as the Nyman-Beurling criterion for the Riemann Hypothesis.
\\\\
Using a fundamental identity:
\begin{equation}
- \frac{\zeta(s)}{s} = s\int_{0}^{\infty} \rho\left(\frac{1}{t}\right) t^{s-1} dt
\end{equation}
where $\rho(x)$ denotes the fractional part of x. Nyman and Beurling formulated the Riemann Hypothesis in terms of the property of the function $\rho$. If $C$ denotes the linear manifold of functions,
$$f(x) = \sum_{k=1}^{n} c_k \rho\left(\frac{\theta_k}{x}\right) 0 < \theta_k \leq 1$$ 
where $c_k$ are constants such that $\sum_{k=1}^{n} c_k \theta_k = 0$. Beurling showed that :\\\\
\textit{The Riemann zeta-function is free from zeros in the half-plane $\sigma > 1/p$ where $1 < p < \infty$ if and only if C is dense in the space $L^p(0,1)$}
\\\\
which is a really wonderful statement that links two separate branches of mathematics in a very elegant manner. This beautiful result prompted further study of Nyman-Beurling's approach. B\'aez-Duarte was among those whose work on this criterion deserves mention. 
\\\\
\noindent
The approach I take in this paper is influenced by Beurling's 1955 paper: \textit{A closure problem related to the Riemann zeta function} and B\'aez-Duarte's 2001 paper: \textit{New versions of the Nyman-Beurling criterion for the Riemann Hypothesis}. In this paper I would be studying the function $\nu$ corresponding to the Dirichlet eta function $\eta$ just like $\rho$ corresponds with the $\zeta$ function. 
\\\\
\section{Dirichlet eta function $\eta(s)$ and $\nu(x)$}

\textbf{Theorem 2.1:}
\textit{For all $\Re(s) > 0$ and $ \nu(t) = \left\{ t/2 \right\} + 1/2 -  \left\{ t/2 + 1/2 \right\} $ we have
\begin{equation}
\eta(s) = s\int_{1}^{\infty} \nu(t) t^{-s-1} dt
\end{equation}
}
\textbf{Proof:}
The Dirichlet eta function is defined as, for all $\Re(s) > 0$,
$$\eta(s) = \sum_{k = 1}^{\infty} \frac{1}{(2k - 1)^s} - \frac{1}{(2k)^s}$$
Now we have for all $t \in \mathbb{R_+}$ $\nu(t) = 0$ or $1$. It is not hard to see that $\nu(t) = 0$ whenever $t \in [2k,2k+1)$ and $1$ whenever $t \in [2k-1,2k)$ for all positive integers $k$. Hence, we can write the integral (2) as
$$s\int_{1}^{\infty} \nu(t) t^{-s-1} dt = \sum_{k=1}^{\infty} \int_{2k-1}^{2k} s t^{-s-1} dt$$
Giving us the sum
$$\sum_{k = 1}^{\infty} \frac{1}{(2k - 1)^s} - \frac{1}{(2k)^s}$$
which is nothing but $\eta(s)$. Since we already know that this sum converges for $\Re(s) > 0$, we get our result.
\\\\
Alternatively, we can write equation (2) as
\begin{equation}
\frac{\eta(s)}{s} = \int_{0}^{1} \nu\left(\frac{1}{t}\right) t^{s-1} dt
\end{equation}
\section{Lemmas and Theorems}
\noindent
\textbf{Note:} Before we proceed further, we set $f(x) = \sum_{k=1}^{\infty} \mu(k) \nu(x/k)$. We shall denote it simply by $f$. In other cases, $f(x) = \sum_{k=1}^{m} \mu(k) \nu(x/k)$ will be denoted by $f_m$. Here, $\mu$ is the usual Moebius-mu function. Also, $\eta(s) = (1 - 2^{1-s}) \zeta(s)$ for $\Re(s) > 0$. 
\\\\
Again, let us keep in mind the following properties of $\nu(x)$ for $x \in \mathbb{R}$:
\\
1. $\nu(x)$ takes only the values 0 and 1.\\
2. $\nu(x)$ is constant in any interval $[n,n+1)$ where $n \in \mathbb{N}$\\
3. $\nu(x/y) = 0$ if $x < y$ for $x,y \in \mathbb{R}$ and $x > 0$ and $y > 0$
\\\\
\noindent
\textbf{Lemma 3.1:} \textit{For $0 < \theta \leq 1$ and $\Re(s) > 0$ we have $\int_{0}^{1} \nu\left(\frac{\theta}{x}\right) x^{s-1} dx = \theta^s \frac{\eta(s)}{s}$}
\\\\
\noindent
\textbf{Proof:}
We have for $0 < \theta \leq 1$
$$ \int_{0}^{\theta} \nu\left(\frac{\theta}{x}\right) x^{s-1} dx = \frac{ \theta^s \eta(s)}{s}$$
If $\theta = 1$ we have our result. Hence suppose that $0 < \theta < 1$.
Since $\nu(\theta/x) = 0$ whenever $ \theta < x$ therefore $$\int_{\theta}^{1} \nu\left(\frac{\theta}{x}\right) x^{s-1} dx = 0 $$
Summing up the above two integrals we get the following equation
\begin{equation} 
\int_{0}^{1} \nu\left(\frac{\theta}{x}\right) x^{s-1} dx = \frac{ \theta^s \eta(s)}{s}
\end{equation}
\\\\
\noindent
\textbf{Lemma 3.2:} For all $x \geq 2$ , $f(x)= -1$
\\\\
\noindent
\textbf{Proof:}
For $\Re(s) > 0$ we get from Lemma 3.1
\begin{equation} 
\int_{0}^{1} f\left(\frac{1}{x}\right) x^{s-1} dx = \frac{\eta(s)}{s}\sum_{k=1}^{\infty} \frac{\mu(k)}{k^s}
\end{equation}
Fix $s$ with $\Re(s) > 1$. For $\Re(s) > 1$ we know that $\sum_{k=1}^{\infty} \mu(k)/k^s = 1/\zeta(s)$.
\\
Hence,
\begin{equation*} 
\int_{0}^{1} f\left(\frac{1}{x}\right) x^{s-1} dx = \frac{1-2^{1-s}}{s}
\end{equation*}
Giving us,
\begin{equation*} 
\int_{0}^{1} \left( 1 + f\left(\frac{1}{x}\right)\right) x^{s-1} dx = \frac{1}{s} +  \frac{1-2^{1-s}}{s}
\end{equation*}

\begin{equation*} 
\int_{0}^{\frac{1}{2}} \left( 1 + f\left(\frac{1}{x}\right) \right) x^{s-1} dx + \int_{\frac{1}{2}}^{1} \left( 1 + f\left(\frac{1}{x}\right) \right) x^{s-1} dx = \frac{1}{s} + \frac{1-2^{1-s}}{s}
\end{equation*}
Now, since $f(x) = 1$ whenever $x \in [1,2)$ giving us for all $\Re(s)  > 1$
Hence, we get
\begin{equation} 
\int_{0}^{\frac{1}{2}} \left( 1 + f\left(\frac{1}{x}\right) \right) x^{s-1} dx + \int_{\frac{1}{2}}^{1}2 x^{s-1} dx = \frac{1}{s} + \frac{1-2^{1-s}}{s}
\end{equation}
ie.,
\begin{equation} 
\int_{0}^{\frac{1}{2}} \left( 1 + f\left(\frac{1}{x}\right) \right) x^{s-1} dx = 0
\end{equation}
Then we have from equation (7) that  for all $\Re(s) > 1$
\begin{equation*}
\int_{\frac{1}{3}}^{\frac{1}{2}} \left( 1 + f\left(\frac{1}{x}\right) \right) x^{s-1} dx + \int_{\frac{1}{4}}^{\frac{1}{3}} \left( 1 + f\left(\frac{1}{x}\right) \right) x^{s-1} dx + ... = 0
\end{equation*}
Since, $f(r)$ is always constant in any given $[r,r+1)$ where $r \in \mathbb{N}$, we get
$$(1 + f(2)) \left(\frac{1}{2^s} - \frac{1}{3^s}\right) + (1 + f(3)) \left(\frac{1}{3^s} - \frac{1}{4^s}\right) + ... = 0$$
It follows from the uniqueness property of a Dirichlet series that each of the terms $1 + f(r) = 0$. Therefore for all $x \geq 2$ , $f(x) = -1$.
\\\\
\noindent
\subsection{The Criterion, two proofs}
\textbf{Theorem 3.1:} \textit{The Riemann Hypothesis is true if
\begin{equation}
\lim_{n\to\infty} \int_{n}^{\infty} \left|1 + \sum_{k=1}^{n} \mu(k) \nu(x/k)\right|^2 x^{-2} dx = 0
\end{equation}
}
\textbf{Proof 1:} 
We get from (5) for $\Re(s) > 0$
\begin{equation*}
\int_{0}^{1} f\left(\frac{1}{x}\right) x^{s-1} dx = \frac{\eta(s)}{s}\sum_{k=1}^{\infty} \frac{\mu(k)}{k^s}
\end{equation*}
Giving us,
\begin{equation*}
\int_{0}^{\frac{1}{2}} f\left(\frac{1}{x}\right) x^{s-1} dx + \int_{\frac{1}{2}}^{1} f\left(\frac{1}{x}\right) x^{s-1} dx = \frac{\eta(s)}{s}\sum_{k=1}^{\infty} \frac{\mu(k)}{k^s}
\end{equation*}
Since $f(x) = 1$ whenever $x \in [1,2)$, we get
\begin{equation*}
\int_{0}^{\frac{1}{2}} f\left(\frac{1}{x}\right) x^{s-1} dx + \int_{\frac{1}{2}}^{1} x^{s-1} dx = \frac{\eta(s)}{s}\sum_{k=1}^{\infty} \frac{\mu(k)}{k^s}
\end{equation*}
Hence for $\Re(s) > 0$
\begin{equation} 
\int_{0}^{\frac{1}{2}} \left(1 + f\left(\frac{1}{x}\right)\right) x^{s-1} dx = \frac{\eta(s)}{s}\sum_{k=1}^{\infty} \frac{\mu(k)}{k^s} - \frac{1 - 2^{1-s}}{s}
\end{equation}
Define $f^*(x) = f(1/x)$ for $0 < x \leq 1/2$. Using Holders inequality we have for $L^2(0,\frac{1}{2})$
\begin{equation}
\left| \frac{\eta(s)}{s}\sum_{k=1}^{\infty} \frac{\mu(k)}{k^s} - \frac{1 - 2^{1-s}}{s} \right| \leq || 1 + f^*||_2 . || x^{s-1} ||_2
\end{equation}
\\\\
If $\Re(s) > \frac{1}{2}$, we get for $s = \sigma + it$
$$|| x^{s-1} ||_2 = \left( \int_{0}^{1/2} |x^{s-1}|^2 dx \right)^{1/2} = \frac{1}{2^{\sigma - \frac{1}{2}} \sqrt{2\sigma - 1}}$$
Now, in equation (10) if $|| 1 + f ^*||_2 = 0$ then the Riemann Hypothesis follows, because it implies $\sum_{n=1}^{\infty} {\mu(n)/n^s}$ converges in the half-plane $\Re(s) > \frac{1}{2}$ and where it exactly equals ${1/\zeta(s)}$. To show $|| 1 + f^* ||_2 = 0$ in $L^2(0,\frac{1}{2})$ we have to prove that
$$\lim_{n\to\infty} \left( \int_{0}^{1/2} \left|1 + \sum_{k=1}^{n} \mu(k) \nu(1/kx)\right|^2 dx \right)^{1/2} = 0$$
But since by Lemma 3.2, $f(x) = \sum_{k=1}^{n} \mu(k) \nu(x/k) = -1$ for $2 \leq x \leq n$, hence we have
\begin{equation}
\lim_{n\to\infty} \int_{0}^{1/2} \left|1 + \sum_{k=1}^{n} \mu(k) \nu(1/kx)\right|^2 dx  = \lim_{n\to\infty} \int_{n}^{\infty} \left|1 + \sum_{k=1}^{n} \mu(k) \nu(x/k)\right|^2 x^{-2} dx
\end{equation}
which completes our proof for the criterion.
\\\\
\noindent
\textbf{Proof 2:}  In this second proof we show that $\sum_{k=1}^\infty {\mu(k)}{/k^s}$ forms a Cauchy sequence for $\Re(s) > \frac{1}{2}$ if 
\begin{equation}
\lim_{n\to\infty} \left( \int_{n}^{\infty} \left|1 + \sum_{k=1}^{n} \mu(k) \nu(x/k)\right|^2 x^{-2} dx \right)^{1/2} = 0
\end{equation}

This simply follows from the following. From lemma 3.1 and using the fact that $\sum_{k=1}^{n} \mu(k) \nu(x/k) = 1$ for $n > 2$ and $x \in [1,2)$ we get for $n > m > 2$
\begin{equation}
\frac{\eta(s)}{s} \sum_{k=m+1}^{n} \frac{\mu(k)}{k^s} = \int_0^{\frac{1}{2}} \left( 1 + \sum_{k=m+1}^n \mu(k)\nu\left( \frac{1}{kx} \right) \right) x^{s-1}dx
\end{equation}
If we let $n \to \infty$ then 
\begin{equation}
\frac{\eta(s)}{s} \sum_{k=m+1}^{\infty} \frac{\mu(k)}{k^s} = \int_0^{\frac{1}{2}} \left( 1 + \sum_{k=m+1}^\infty \mu(k)\nu\left( \frac{1}{kx} \right) \right) x^{s-1}dx
\end{equation}
But using the result of Lemma 3.2 in the RHS of (14) equals
\begin{equation}
\begin{split}
& \int_0^{\frac{1}{2}} \left( 1 + \sum_{k=m+1}^\infty \mu(k)\nu\left( \frac{1}{kx} \right) \right) x^{s-1}dx\\
& = \int_0^{\frac{1}{2}} \left( 1 + \sum_{k=1}^\infty \mu(k)\nu\left( \frac{1}{kx} \right) \right) x^{s-1}dx -  \int_0^{\frac{1}{2}} \left( 1 + \sum_{k=1}^m \mu(k)\nu\left( \frac{1}{kx} \right) \right) x^{s-1}dx\\
& = -  \int_0^{\frac{1}{2}} \left( 1 + \sum_{k=1}^m \mu(k)\nu\left( \frac{1}{kx} \right) \right) x^{s-1}dx
\end{split}
\end{equation}
Using a similar treatment as in proof 1, we get  by applying Holder's inequality for $L^2(0,\frac{1}{2})$ to equation (14) and replace the RHS by (15) to get,
\begin{equation}
\left|\frac{\eta(s)}{s} \sum_{k=m+1}^{\infty} \frac{\mu(k)}{k^s} \right| \leq || 1 + f^*_m ||_2. || x^{s-1} ||_2
\end{equation}
Where $$f_m\left(\frac{1}{x}\right) =  \sum_{k=1}^m \mu(k)\nu\left( \frac{1}{kx} \right) \quad \text{and} \quad f^*_m(x) = f_m\left(\frac{1}{x}\right)$$
Here, if $|| 1 + f^*_m ||_2 \to 0$ in $L^2(0,1/2)$ as $m \to \infty$ then  we can say that $\sum_{k=1}^{n} {\mu(k)}/{k^s}$ forms a Cauchy sequence and hence converges for $\Re(s) > 1/2$ and this implies the Riemann Hypothesis. 
The criterion is then derived in a similar fashion as in the first proof.

\end{document}